\documentclass[journal]{IEEEtran}
\usepackage[T1]{fontenc}
\usepackage{lmodern}
\usepackage{amsmath}
\usepackage{booktabs}
\usepackage[inline]{enumitem}
\usepackage{graphbox}
\usepackage{tabularx}
\usepackage{cleveref}
\usepackage[noend]{algpseudocode}
\usepackage{algorithm}
\usepackage{multicol}
\usepackage{comment}
\usepackage{multirow}
\usepackage{todonotes}
\usepackage{graphicx}
\graphicspath{ {./Images/} }
\usepackage{caption}
\usepackage{subcaption}
\usepackage{url} 
\begin{document}
\title{Comparing Optimization Models for Radiotherapy Scheduling}
\author{
    \IEEEauthorblockN{
    {Chiara Camilla}~{Rambaldi Migliore}\IEEEauthorrefmark{1,2}, 
    {David}~{Stanicel}\IEEEauthorrefmark{1}, 
    {Nysret}~{Musliu}\IEEEauthorrefmark{3}, 
    {Giovanni}~{Iacca}\IEEEauthorrefmark{1} 
    and {Marco}~{Roveri}\IEEEauthorrefmark{1}
    }
    \IEEEauthorblockA{\IEEEauthorrefmark{1}University of Trento, Italy}
    \IEEEauthorblockA{\IEEEauthorrefmark{2}University of Pisa, Italy}
    \IEEEauthorblockA{\IEEEauthorrefmark{3}Technische Universität Wien, Austria}
}

\markboth{C.C. Rambaldi Migliore et al.}{Comparing Optimization Models for Radiotherapy Scheduling}

\maketitle

\begin{abstract}
The Radiotherapy Scheduling Problem (RTSP) involves determining an optimal schedule for patients undergoing radiation treatments, a task that has a massive impact on clinical outcomes given the central role of radiotherapy in cancer care. The daily batch approach--which consists of scheduling all the newly arrived patients together at the end of each day--modelled with Integer Linear Programming, is currently one of the most effective methods for the RTSP. However, this kind of formulation requires substantial computational resources in terms of time and memory. Here, we address these limitations by developing two novel greedy heuristics (named RTSP First Fit and RTSP Best Fit) and use them as constructive heuristics for a Simulated Annealing (SA) approach to optimize the scheduling. The proposed methods--the heuristics alone and their combination with SA--are evaluated on a publicly available dataset against an integer linear program formulation solved with two different state-of-the-art exact solvers. Evaluation metrics include six scheduling objectives capturing patient waiting times, preference satisfaction, and changes in linear accelerator assignment (aggregated in four different weight configurations), solving time, and memory consumption. The results show that the novel heuristics achieve solutions close to those of exact methods, while dramatically reducing runtime and memory usage; furthermore, when combined with SA, they further improve the solution quality while maintaining low runtime and memory usage.
\end{abstract}

% Note that keywords are not normally used for peer review papers.

\section{Introduction}

Radiotherapy is currently one of the most effective cancer treatments \cite{Roberts2023RadiationTF}. It aims to deliver precise doses of radiation to tumours while preserving healthy tissue. In addition, early initiation of treatment is crucial for achieving optimal results. Therefore, effective scheduling of therapies is essential to balance patient needs, resource constraints, and treatment efficacy.

The research community addressed the Radiotherapy Scheduling Problem (RTSP) \cite{vieira2016operations} by introducing various formulations and algorithms, ranging from exact Operational Research (OR) methods \cite{pham2022two,frimodig2023comparing,frimodig2023column} to metaheuristics \cite{cares2013ls2r,vogl2019scheduling,vogl2017multi,braune2022stochastic}.
Existing solutions can be roughly divided into three main categories, namely: (i) \emph{offline scheduling} approaches, which search for optimal appointments for all patients at the same time; (ii) \emph{batch scheduling} approaches, where newly arrived patients are grouped (daily or weekly) and then scheduled, providing the scheduling procedure information about only a small group of patients; and (iii) \emph{online scheduling} approaches, where the goal is to schedule one patient at a time, without knowledge of future patients.
Besides these formulations, only a few previous works have attempted to map the RTSP to existing combinatorial problems, such as the Bin Packing Problem (BPP). Furthermore, we are not aware of any previous attempt leveraging heuristics to rapidly generate good solutions for the RTSP, and then feeding them to a metaheuristic, such as Simulated Annealing (SA), to achieve good-quality results with limited resources. Here, we address these gaps by making the following contributions:
\begin{itemize}[leftmargin=*,nosep]
    % \item First, inspired by \cite{migliore2024bin,migliore2025addressing}, we adopt a 1D BPP formulation as a basis to model the RTSP, extending and reinterpreting their work to incorporate new problem variables. 
    % %The 1D BPP is an NP-hard combinatorial optimisation problem \cite{garey1979computers} that involves packing items of various sizes into a minimum number of identical-size bins. 
    % By mapping radiotherapy treatments to items and a machine (i.e., a linear accelerator (LINAC)) to a bin, we directly apply the 1D BPP concepts to the RTSP.
    \item We adopt the 1D BPP RTSP formulation from \cite{migliore2024bin} to develop two novel heuristics, inspired by the First Fit (FF) and Best Fit (BF) heuristics \cite{johnson1974worst}.
    %we develop new adaptations of the First Fit and Best Fit \cite{johnson1974worst} heuristics for the 1D BPP to address the RTSP.
    \item We develop an SA method to further improve the heuristic solutions and evaluate three distinct variants.
    %We develop three different SA implementations to optimize the solutions given by the two heuristics.
    \item We perform a thorough experimental evaluation of the proposed methods on a public dataset simulating patient arrivals, evaluating performance across (i) four different combinations of scheduling objectives (capturing patient waiting times, preference satisfaction, and changes in machine assignments), (ii) solving time, and (iii) memory consumption. We compare heuristics alone, and their combination with SA, against the Integer Linear Programming (ILP) formulation proposed in \cite{frimodig2023comparing}
    %--whose details are provided in the Supplementary Material--
    solved using both the CPLEX commercial solver and the OR-Tools CP-SAT open-source solver with different problem configurations. This comparison provides practical indications about the most suitable approach to choose depending on the available resources. 
    %computational resources, including the cost of purchasing a license for a commercial solver.
    %\footnote{Despite the several exchange of emails with the authors of \cite{frimodig2023comparing} for clarifications, the code was not released or made available to us for comparison.} 
\end{itemize}
The results show that the proposed heuristics can provide the SA with good solutions from which to start the search toward the optimal solution. Moreover, both the heuristics alone and the SA implementations perform remarkably better in terms of runtime and memory usage compared to the ILP baseline formulation. This improvement is achieved while preserving a good solution quality in comparison to the exact solvers. This is an important finding, for a number of reasons, as: (1) physicians may not have dedicated servers on which they can run exact resource-intensive methods; (2) OR solvers are not always open-source, and (3) hospitals may not have enough funds for software licences.

%%%%%%%%%%%%%%%%%%%%%%%%%%%%%%%%%%%%%%%%%%%%%%%%%%%%%%%%%%%%%%%%%%%%%%%%

\section{Related Works}
\label{sec:related_works}

%Scheduling in the healthcare field is a broader topic that encompasses several specific problems, like e.g., bed capacity \cite{zhang2012simulation, holm2013improving}, operating room capacity \cite{wang2009sizing}, Transcranial Magnetic Stimulation scheduling \cite{squires2022novel}, all addressed with several techniques (e.g., exact methods, simulation-based, genetic algorithms). 

Scheduling in healthcare is a broad topic that encompasses several specific problems. For the purpose of this paper, we focus on scheduling radiotherapy appointments, categorising approaches into 
\begin{enumerate*}[label=(\arabic*)]
\item offline approaches, 
\item batch approaches, and 
\item online approaches.
\end{enumerate*}
For each category, we first analyse the exact operations research methods, followed by the approaches based on metaheuristics.

\paragraph{Appointment scheduling with offline approaches} 
Several studies have addressed offline radiotherapy appointment scheduling, often using OR methods and metaheuristics while bypassing pre-treatment phases. 
Early work by Conforti et al. \cite{conforti2008optimization} formulated the RTSP as an ILP problem first to maximise the number of patients starting treatment and then to minimise patient waiting times \cite{conforti2010non} and patient availability \cite{conforti2011optimal}, testing their approaches on both toy and real-world data.
Vogl et al. \cite{vogl2017multi,vogl2019scheduling} studied metaheuristics—Genetic Algorithm (GA), Iterated Local Search (ILS), and a hybrid of the two—for an ion beam facility, modelling the task as a modified job shop scheduling problem with custom constraints to minimise bottleneck resource time and violations of daily \emph{time window} constraints (e.g., early/late morning, early/late afternoon).
%Metaheuristic approaches, namely a Genetic Algorithm (GA), Iterated Local Search (ILS), and a hybrid of the two, have been studied by Vogl et al. \cite{vogl2017multi,vogl2019scheduling}, who approached the task at an ion beam facility by modelling it as a modified job shop scheduling problem with custom constraints, aiming to minimise bottleneck resource time and violations of \emph{time window} (e.g., early/late morning and early afternoon/late afternoon) defining the daily schedule.

\paragraph{Appointment scheduling with batch approaches} 
%Several authors adopted batch scheduling strategies, where patient appointments are planned over a multi-day or weekly horizon. Existing works on batch scheduling propose ILP, Mixed ILP (MILP), Column Generation (CG) approaches, or combinations thereof, and approaches based on metaheuristics.
Numerous researchers have employed batch scheduling strategies, where patient appointments are arranged over a multi-day or weekly time frame. Previous studies on batch scheduling suggest the utilisation of ILP, Mixed ILP (MILP), Column Generation (CG) methods, or combinations of these, along with metaheuristics.

Jacquemin et al. \cite{jacquemin2010towards} introduced an ILP model to schedule patients weekly (over a 15-week period), using a small dataset with two identical machines. They considered patient availability and continuity of care,% and later refined the model by incorporating,
treatment patterns, optimising the use of the radiotherapy room, treating more patients, and reducing waiting times \cite{jacquemin2011pattern}. \cite{burke2011integer} used an ILP model daily with various optimisation criteria, determining that scheduling within 7 days post-release was optimal. \cite{migliore2024bin,migliore2025addressing} examined an ILP model cast as a BPP on small synthetic instances and compared it with a GA using real data. %These BPP-based approaches inspired this work.

Vieira et al. \cite{vieira2020radiotherapy} introduced a MILP formulation with window preferences, employing heuristics for machine pre-assignment, and then using MILP for weekly scheduling. \cite{vieira2021radiotherapy} improved this by incorporating patient time preferences, reducing machine switches, and aligning more closely with patient preferences than manual scheduling. \cite{boonmee2021decision} introduced constraints involving doctor availability and treatment coordination, proving feasibility with data from the Thai Cancer Center. \cite{emsamrit2024mixed} refined the model by adding surgical recovery times and optimising capacity constraints, improving computational efficiency. Hybrid approaches that combine ILP/MILP and Constraint Programming (CP) were explored, such as the two-phase method proposed in \cite{pham2022two}, including (i) session assignment using ILP and (ii) scheduling with MILP and CP. 
%verified on real data, favouring batch over online scheduling.

Saure et al. \cite{saure2012dynamic} initially explored the CG approach by modelling the RTSP as a discounted infinite-horizon Markov Decision Process, converting it into an LP, and implementing CG daily to manage appointments and overtime. Frimodig et al. \cite{frimodig2023comparing} evaluated ILP, CG-ILP, and CP to reduce waiting times and align with patient preferences. Later, \cite{frimodig2023column} improved the CG model to deal with machine unavailability, establishing the most comprehensive public radiotherapy dataset to date. Their ILP serves as the baseline for our paper.

A metaheuristic-based approach was introduced in \cite{cares2013ls2r}, proposing an effective Local Search (LS) method for the RTSP. Maschler et al. \cite{maschler2016particle} addressed the particle therapy scheduling problem, initially using MILP but finding it intractable. They then developed an algorithm with a forward-looking mechanism and proposed a Greedy Randomised Adaptive Search Procedure (GRASP) along with an iterated greedy metaheuristic for daily patient scheduling.
Braune et al. \cite{braune2022stochastic} employed a GA-based approach and Monte Carlo simulations to manage appointment durations with variable treatment times, with a focus on uncertainty in patient preparation and exit times for daily scheduling.

\paragraph{Appointment scheduling with online approaches} 
To our knowledge, the only work that proposed an online scheduling approach for RTSP, in which patients are individually scheduled, was done in \cite{legrain2015online}, incorporating patient arrival uncertainty into a hybrid stochastic online model to enhance scheduling.

\paragraph{Novelty of this work} 
To the best of our knowledge, this is the first work to tackle the 1D BPP formulation of the RTSP using newly designed greedy heuristics inspired by online 1D BPP methods, and the first to explore how these heuristics can be further sharpened via SA. Moreover, it is the first study to conduct a systematic, head-to-head comparison of these heuristic and metaheuristic approaches against exact methods on a common dataset.
%To the best of our knowledge, no previous research has addressed the 1D BPP formulation of the RTSP using newly devised greedy heuristics inspired by online 1D BPP methods, nor has it examined how these heuristics can be refined through SA. Furthermore, this is the first study to systematically compare these heuristic and metaheuristic strategies with exact methods on the same dataset.

%%%%%%%%%%%%%%%%%%%%%%%%%%%%%%%%%%%%%%%%%%%%%%%%%%%%%%%%%%%%%%%%%%%%%%%%

\section{Problem formulation}
\label{sec:prob_form}

The RTSP consists of scheduling appointments for all patients' \textbf{fractions} (single radiotherapy sessions delivering part of the total dose), i.e., assigning each fraction a \textbf{day} and a \textbf{time window} (a predefined part of the day), on each specific \textbf{machine}, subject to some constraints. The formulation and solution quality metric presented here are taken from \cite{frimodig2023comparing}.%, since we use the same dataset and ILP specification for experiments and comparison.

%\subsection{Problem Input}
\subsubsection{Problem Input.} The problem has three main inputs: a list of \textbf{patients}, a \textbf{time horizon}, and a list of \textbf{machines}.
Each patient has a given number of fractions to schedule and an associated treatment protocol, which defines the priority of the treatment and the preferred and allowed machines on which the fraction can be performed. Each patient also has a minimum start day--which from now on will be called the ``arrival day''--and a target day, i.e., the maximum day after which it is considered late to start the treatment. Patients may or may not have a preference for the time window of the day in which to be scheduled (for example, early/late morning or afternoon).
The time horizon is defined by the number of days considered for the scheduling.
Each machine has a residual time capacity for each time window of each day of the time horizon. Some machines may be partially or fully beam-matched, meaning that the same radiotherapy plan can be safely administered by any of them. Machines are considered fully matched when they are located in the same building; if they are in different buildings, they are partially matched.

%\subsection{Solution}
\subsubsection{Solution.}
A solution to the RTSP consists of the complete appointment calendar for each patient. 
%A solution is considered feasible if it satisfies the following criteria:
A feasible solution satisfies the following criteria:
\begin{enumerate}[label=\arabic*),nosep]
    \item \textbf{Singularity}: There cannot be more than one fraction per day; since a day is divided into several time windows, each fraction must be scheduled only in one time window of a machine for each day.
    \item \textbf{Consecutiveness}: For each patient, fractions must be scheduled consecutively on weekdays.
    \item \textbf{Availability}: A fraction can be scheduled on a machine in a time window only if there is enough time remaining for that machine in that time window.
    \item \textbf{Starting}: A patient can only start treatment on protocol-allowed days, not before their arrival, and no later than the last day of the time horizon minus the total number of fractions to be scheduled.
   % A patient cannot start the treatment on days that are not allowed by their protocol, before their arrival day, and after a maximum day calculated as the last day in the time horizon minus the total number of fractions to be scheduled for the patient.
    \item \textbf{Specificity}: The fractions of each patient can be performed only on the allowed machines specific to their treatments.
    \item \textbf{Protocol precedence}: For two patients sharing the same treatment protocol, the one with the earlier target day must start no later than the other.
\end{enumerate}

%\subsection{Evaluation Metrics}
\subsubsection{Evaluation Metrics.}
%The evaluation metrics used for the comparative analysis of the experiments are solution quality, solving time, and memory usage. The solution quality is measured as a weighted sum of six objectives--plus an offset of 1, as in \cite{frimodig2023comparing}:
The comparative analysis uses three evaluation metrics: solution quality, solving time, and memory usage. Solution quality is computed as a weighted sum of six objectives, plus an offset of 1, as in \cite{frimodig2023comparing}:
\begin{equation}%\small
\label{eq:objective_sum}
    \scalebox{0.85}{$\text{Obj. Sum} = 1 + \alpha_1 f_1 + \alpha_2 f_2 + \alpha_3 f_3 + \alpha_4 f_4 + \alpha_5 f_5 + \alpha_6 f_6$}
\end{equation}

\noindent Each objective is described as follows:
\begin{itemize}[leftmargin=*,nosep]
    \item $f_1$: the weighted sum of days patients are waiting before starting the treatment (the weights are 10, 3, 1, respectively, for priority 1, 2, and 3, where the lower value for priority means higher precedence);
    \item $f_2$: the weighted sum of days beyond the deadline for patients to start their treatment (the weights are the same as for $f_1$);
    \item $f_3$: the number of times patients switch time windows between two fractions;
    \item $f_4$: the sum of the deviations from the time window preferred by the patients over all fractions;
    %the number of times patients' fractions are scheduled on a non-preferred machine indicated by the treatment protocol;
    \item $f_5$: the number of patients' fractions scheduled on a non-preferred machine as per the treatment protocol;
    \item $f_6$: the sum of switches between machines that are partially beam-matched, i.e., machines whose treatment plans are interchangeable (they can be thought of as clones), but located in different buildings.
\end{itemize}

We denote by \#n : $(\alpha_1,\alpha_2, \alpha_3,\alpha_4,\alpha_5,\alpha_6)$ a combination of weights. As in \cite{frimodig2023comparing}, we consider, in total, the following four combinations:
{
\setlength{\abovedisplayskip}{0pt}
\setlength{\belowdisplayskip}{0pt}
\setlength{\abovedisplayshortskip}{0pt}
\setlength{\belowdisplayshortskip}{0pt}
\begin{align*}
&\#1 : (50,100,1,0,10,10) & &\#2 : (50,100,1,1,0,0)\\
&\#3 : (100,0,1,0,10,0)   & &\#4 : (100,0,1,5,10,10).
\end{align*}
}

%\noindent Combination \#1 can fit better those hospitals that do not consider patients' preferences for treatment times, while Combination \#2 is suitable for cancer centres that do not have multiple buildings and do not state a preferred machine for treatments. The combination \#3 may suit hospitals that do not have target days, do not consider preferences, and do not have preferred machines. Lastly, Combination \#4 considers all objectives but the deadlines, therefore fitting cancer centers that consider everything but the deadlines.

\noindent Combination \#1 suits hospitals that do not consider patients' preferred treatment times. Combination \#2 suits cancer centres without multiple buildings and without preferred machines. Combination \#3 suits hospitals with no target days, no patient preferences, and no preferred machines. Combination \#4 considers all objectives except deadlines and thus fits cancer centres that consider everything but deadlines.

%%%%%%%%%%%%%%%%%%%%%%%%%%%%%%%%%%%%%%%%%%%%%%%%%%%%%%%%%%%%%%%%%%%%%%%%

\section{Methodology}
This section presents the proposed methodology used to tackle the RTSP. In particular, we briefly recall the ILP model, then describe the two proposed heuristics and the SA, along with its different implementations.

\subsection{Integer Linear Programming Model}
In order to compare our approach with the state-of-the-art formulation from \cite{frimodig2023comparing}, we used their ILP formulation and implemented the model with both the Google OR-Tools framework and the CPLEX Python API.
The problem is formulated with three principal pseudo-Boolean variables for each patient $p$, namely: (1) $q_{p,m,d,f}$, which takes value 1 if the fraction $f$ of the patient is scheduled on day $d$ and machine $m$ and 0 otherwise; (2) $x_{p,m,d,w}$, which takes value 1 if the patient has a fraction scheduled for day $d$ on machine $m$ in the time window $w$ and 0 otherwise; (3) $t_{p,m,d,w}$, which can be described as $x_{p,m,d,w}$ but only for the first fraction. Those variables are used to model the constraints and objective function presented in the previous section, along with additional helper variables. For further details on the formulations, see \cite{frimodig2023comparing}.

\subsection{RTSP Heuristics}
\begin{figure}[ht!]
    \centering
    \includegraphics[width=0.95\columnwidth]{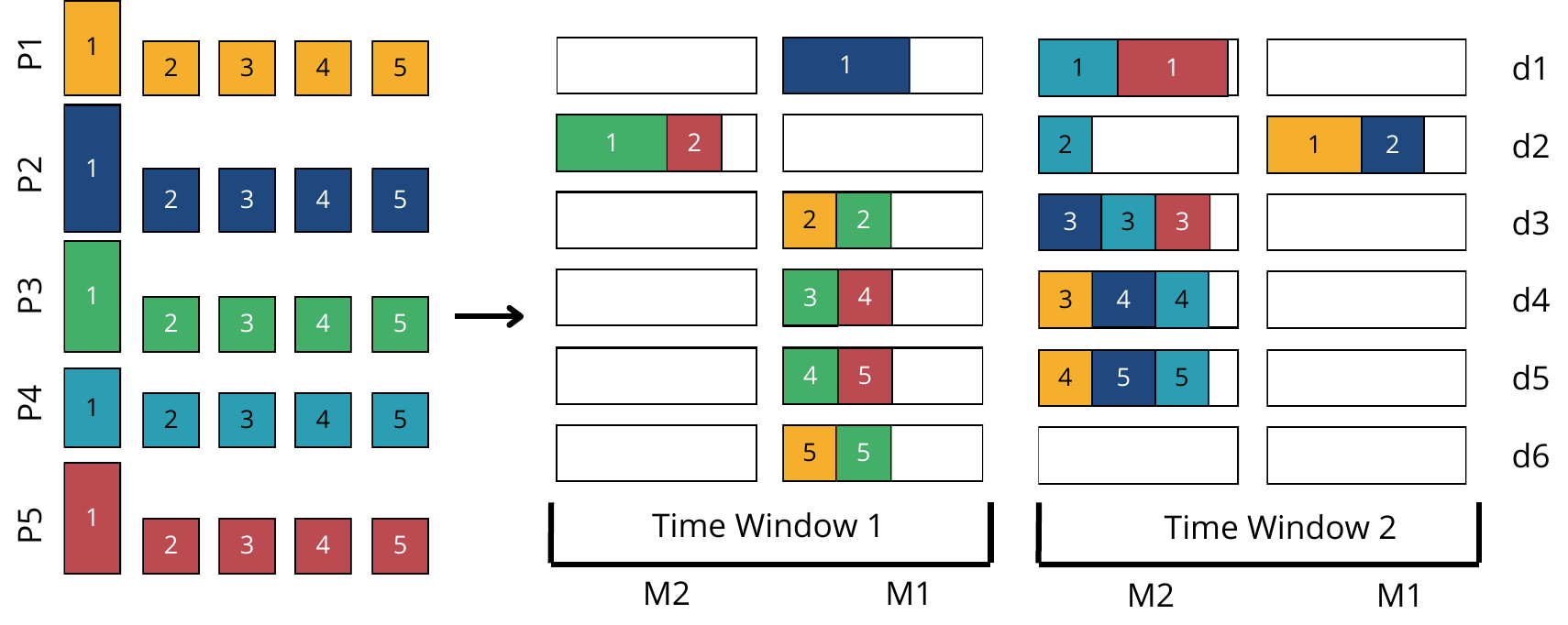}
    \caption{A visual example of the modified 1D BPP. Each colour indicates a patient. Each patient's fraction must appear no more than once in each bin's group.}
    \label{fig:m_bpp}
\end{figure}

Extending the work from \cite{migliore2024bin}, we reinterpret the RTSP as a modified 1D BPP by incorporating time windows. We provide a visual representation of this mapping in \Cref{fig:m_bpp}. 
According to this modified 1D BPP, items correspond to fractions and are grouped by patient, while bins correspond to machines operating during a specific time window, grouped by day.

\begin{comment}
    In accordance with the constraints described in the previous section, we have the following. 
    \begin{itemize}[leftmargin=*,nosep]
        \item Each item in a given group must be placed contiguously in the groups of bins, as the treatments of a given patient must be on consecutive days.
        \item Each item must be placed in at most one bin for each group, as each treatment must occur only once a day.
        \item Each item placed must not exceed the bin capacity, as each temporal window lasts a finite number of minutes.
        \item The first item of each group can be placed only in the allowed groups of bins, as a patient may need to start the treatment on some specific days.
        \item Each item must be placed in the allowed bins, as patients may need specific machine for their treatment.
    \end{itemize} 
\end{comment}

Based on this analogy, we design two new greedy heuristics drawing inspiration from those developed for the 1D BPP \cite{coffman1984approximation}. Both heuristics deal with one group of items at a time, store the items in groups of bins, and have a mechanism to keep track of the consecutiveness in placing the items. Patients (group of items) are ordered by arrival day, treatment priority, and target day (the maximum number of days they can wait before starting treatment). 
%This ordering resembles the reality of a hospital that, at the end of each day, schedules all new patients of the day one by one without using the knowledge of the others (leveraging the online nature of the heuristics).
%The pseudocodes and source codes of both heuristics are available in the Supplementary Material.

\paragraph{RTSP First Fit (FF)} 
The algorithm starts by ordering the patients by arrival day, treatment priority, and target days, so that they can then be selected in an order that resembles the real one used by the physicians (for each day, the patients with higher priority and a closer target day are scheduled first). Once a patient (corresponding to a group of items) is selected, the heuristic tries to find the first $n$ consecutive groups of bins, where $n$ is the number of items in the selected group (patient), so that in each $i$-th group of bins, at least one allowed bin has enough capacity to store the $i$-th item. In the general FF, if $n$ consecutive groups of bins cannot be found, new groups of bins should be added; however, in our case, the algorithm starts with enough groups to store all the items (this number is taken from the dataset instances), so this part is unnecessary. %The pseudocode of the method is provided in \Cref{code:firstfitmod}.

\begin{comment}
\begin{algorithm}[ht!]
\caption{RTSP First Fit Heuristic}
\label{code:firstfitmod}
    \begin{algorithmic}[1]
        \State $PatientsList \gets$ patients sorted by arrival day, treatment priority, and target day
        \State $allDays \gets$ all days in the scheduling domain
        \State $bins \gets$ available space in each day, machine, window
        \For{$p \in PatientsList$} 
            \For{$day \in allDays$}
                \If{fractions of $p$ consec. in $bins[day]$}
                    \State Store fractions of $p$ in
                    \Statex \ \ \ \ \ \ \  \ \ \ \ \ \ \ \ \ \ \ available $bin[day][machine][window]$
                    \Statex \ \ \ \ \ \ \  \ \ \ \ \ \ \ \ \ \ \ starting from $day$ till 
                    \Statex \ \ \ \ \ \ \  \ \ \ \ \ \ \ \ \ \ \ $day+p.fractionNumber$
                    \State \textbf{break} %Break the days loop
                \EndIf
            \EndFor
        \EndFor
    \end{algorithmic}
\end{algorithm}
\end{comment}

\paragraph{RTSP Best Fit (BF)}
As in RTSP FF, this algorithm first orders patients by arrival day, treatment priority, and target days. Then, it attempts to store the selected group of items (i.e. patient's fractions) in the $n$ consecutive groups of bins, so that in each $i$-th group, the chosen bin is the one with the least remaining capacity and that can accommodate the $i$-th item. Similar to the RTSP FF heuristic, this algorithm starts with a sufficient number of groups (this fixed number is taken from the dataset instances) to store all items by design. This entails that adding new groups of bins is not necessary. 
%The pseudocode of the method is provided in \Cref{code:bestfitmod}.

\begin{comment}
\begin{algorithm}[ht!]
\caption{RTSP Best Fit Heuristic}
\label{code:bestfitmod}
    \begin{algorithmic}[1]
        \State $PatientsList \gets$ patients sorted by arrival day, treatment priority, and target day
        \State $allDays \gets$ all days in the scheduling domain
        \State $bins \gets$ available space in each day, machine, window
        \For{$p \in PatientsList$} 
            \For{$day \in allDays$}
                \For{$f \in p.fractions$}
                    \If{$f$ fit in $bins[day+f-1]$}
                        \State Calculate remaining capacity 
                        \State after $f$ has been added
                    \EndIf
                    \State Store $f$ in $bins[day][machine][window]$ 
                    \Statex \ \ \ \ \ \ \  \ \ \ \ \ \ \ \ \ \ \ with minimum remaining capacity after \Statex \ \ \ \ \ \ \  \ \ \ \ \ \ \ \ \ \ \ storing $f$
                    \If{no such bin exists}
                        \State \textbf{break} %Break the fractions loop
                    \EndIf
                    %\State If no such bin exists, break the fractions loop
                \EndFor
                \If{all fractions are stored consecutively}
                    \State \textbf{break} %Break the days loop
                \EndIf
                %\State If all fractions have been stored consecutively
                %\State break days loop
            \EndFor
        \EndFor
    \end{algorithmic}
\end{algorithm}
\end{comment}

\paragraph{Complexity analysis}
Indicating with $p$ the number of patients, $d$ the number of days, and $f_{max}$ the maximum number of fractions prescribed to patients among all patients, both FF and BF have a complexity of $\mathcal{O}\left(\max (p \cdot log(p), p \cdot d \cdot f_{max})\right)$. The $p \cdot log(p)$ term corresponds to the complexity of the sorting of \textit{PatientsList}, while $p \cdot d \cdot f_{max}$ is the complexity of the \textit{for} loops. From a wide analysis of the RTSP instances, it results that $d \cdot f_{max}$ is always greater than $p \cdot log(p)$. We can therefore conclude that, for both algorithms, the complexity is $\mathcal{O}\left(p \cdot d \cdot f_{max}\right)$.

%%%%%%%%%%%%%%%%%%%%%%%%%%%%%%%%%%%%%%%%%%%%%%%%%%%%%%%%%%%%%%%%%%%%%%%%

\subsection{RTSP Simulated Annealing (SA)}
The two previously introduced heuristics are further improved by SA, which uses the solution found by either of the two heuristics as a starting point for optimisation. As discussed in \cite{sarhani2023initialization}, SA can benefit from starting from solutions constructed with heuristics to avoid a wasteful initial search that would otherwise be needed starting from random solutions. 
SA was chosen because it can escape local minima by occasionally accepting worsening moves and because it is easy to start its search from a solution found by heuristics. Furthermore, it is widely used and studied in combinatorial problems, where one can include different kinds of domain-specific moves, such as the case of the RTSP.
Three different implementations of the algorithm have been developed: \textbf{plain}, \textbf{daily}, and \textbf{reheat}.
Each algorithm takes as input the initial temperature $t_{start}$, the cooling factor $\alpha$, four parameters serving as the size of the neighbourhood of possible moves, and five parameters for the probability of each move being chosen; the maximum number of iterations $k_{max}$ is fixed. All shift operations' resulting indexes are computed modulo the corresponding domain size. The moves are the same for each implementation:
\begin{itemize}[leftmargin=*,nosep]
    \item \textbf{Shift Time Window (m0)}: the move shifts by $k$ the time window assigned to a randomly chosen patient on $j$ days. The shift amount $k$ is decided as $k \gets randomInt(1, tsm)$, where $tsm$ is the \textit{time window shift} parameter. The number of days $j$ is decided as $j \gets randomInt(1, tdm)$, where $tdm$ is the \textit{time window days} parameter.
    \item \textbf{Shift Machine (m1)}: this move shifts by $m$ the machine assigned to a randomly chosen patient on $j$ days. The shift amount $m$ is decided as $m \gets randomInt(1, msm)$, where $msm$ is the \textit{machine shift} parameter. The number of days $j$ is decided as $j \gets randomInt(1, mdm)$, where $mdm$ is the \textit{machine days} parameter;
    \item \textbf{Swap Time Windows (m2)}: two random patients are chosen for swapping their assigned time windows on the same--randomly picked--day;
    \item \textbf{Swap Machines (m3)}: two random patients are chosen for swapping their assigned machines on the same--randomly picked--day;
    \item \textbf{Shift Start Day (m4)}: a random patient is chosen for shifting--backward if possible, forward otherwise--their first day of treatment; all the consequent fractions are then shifted accordingly.
\end{itemize}

\begin{comment}
\begin{algorithm}[ht!]
\caption{Simulated Annealing}
\label{code:SA_plain}
    \begin{algorithmic}[1]
        \State $t \gets t_{start}$
        \For{$k_{max}$} 
            \State $obj_{last} \gets Solution.obj$ 
            \State $move \gets getRandomMove()$
            \State $Solution \gets doMove()$
            \State $obj_{curr} \gets Solution.obj$
            \If{not $AcceptMove(t, obj_{last}, obj_{curr})$}
                \State $undoMove()$
            \EndIf
            \State $t \gets t \cdot \alpha$
        \EndFor
    \end{algorithmic}
\end{algorithm}
\end{comment}

The acceptance criterion is always the standard one: accept a new solution if it is better than the current one, and accept it with probability $e^{-(obj_{curr}-obj_{last})/t}$ if it is worse. 
The SA \textit{plain} implementation corresponds to the standard SA algorithm, as typically presented in foundational literature \cite{kirkpatrick1983optimization}.
The SA \textit{daily} algorithm follows a day-by-day batching logic: the heuristic (fixed at the beginning of the algorithm--either RTSP FF or RTSP BF) generates a schedule for the patients assigned to a day $d$. This schedule is then optimised through SA, machine occupancies are updated, and the algorithm proceeds to the next day. This process is repeated sequentially until all daily batches of patients have been scheduled.
%
%The \textbf{daily} algorithm follows a day-by-day batching logic: as illustrated in \Cref{code:SA_daily}, the heuristic (fixed at the beginning of the algorithm--either RTSP FF or RTSP BF) generates a schedule for the patients assigned to a day $d$. This schedule is then refined through SA. After optimisation, machine occupancies are updated, and the algorithm proceeds to the next day. This process is repeated sequentially until all daily batches of patients have been scheduled.
%
The SA \textit{reheat} algorithm follows the same logic as the plain one, but the temperature is restored after $k_{max}/10$ iterations, hence 10 times, during the course of the run.

\begin{comment}
\begin{algorithm}[ht!]
\caption{Daily Scheduling SA}
\label{code:SA_daily}
\begin{algorithmic}[1]
\State $machOcc \gets InitialOcc()$
\For{$patients \in PatientsDailyBatches$}
    \State $\mathcal{S}_d \gets HeuristicSchedule(patients, machOcc)$
    \State $\mathcal{S}_d \gets SA_{plain}(\mathcal{S}_d, machOcc, k_{max})$
    \State $machOcc \gets UpdateOcc(\mathcal{S}_d)$
\EndFor
\end{algorithmic}
\end{algorithm}
\end{comment}

%%%%%%%%%%%%%%%%%%%%%%%%%%%%%%%%%%%%%%%%%%%%%%%%%%%%%%%%%%%%%%%%%%%%%%%%

\section{Evaluation}
\label{sec:evaluation}

%Here, we first present the dataset, the experimental setup, and then discuss the results.

\subsection{Dataset and experimental setup}

\paragraph{Dataset} We use a public synthetic dataset generated by \cite{frimodig2023comparing} based on historical real-world data collected from Iridium Netwerk, Belgium’s largest cancer centre. The dataset follows a Poisson distribution and is composed of 80 instances, divided into four setups, which are organised as follows (with $\lambda$ being the average arrival rates and $W$ the time windows per day): 20 instances are generated with $\lambda = 16$ and $W = 2$; 20 with $\lambda = 16$ and $W = 4$; 20 with $\lambda = 18$ and $W = 2$; and 20 with $\lambda = 18$ and $W = 4$.
% \begin{itemize}[leftmargin=*,nosep]
%     \item 20 instances generated with $\lambda = 16$ and $W = 2$;
%     \item 20 instances generated with $\lambda = 16$ and $W = 4$;
%     \item 20 instances generated with $\lambda = 18$ and $W = 2$;
%     \item 20 instances generated with $\lambda = 18$ and $W = 4$.
% \end{itemize}
%Each instance in the dataset is provided as a text file containing a comprehensive set of attributes.

\paragraph{Experimental Setup}
We executed all the experiments on a Linux workstation with 256~GB of RAM, two Intel(R) Xeon(R) Gold 6238R CPUs, each featuring 28 physical cores with a base clock speed of 2.20 GHz.

We implemented the proposed heuristics and SA using Python v3.11.11. The baselines chosen for comparison were the ILP formulation solved with the CPLEX solver (hereafter referred to as ILP CPLEX) and with the Google OR-Tools CP-SAT (hereafter referred to as ILP CP-SAT). We implemented these encodings to the best of our knowledge, as the original code is not publicly available\footnote{We had several exchanges of emails with the authors of \cite{frimodig2023comparing} for clarifications; however, the code was neither released nor made available to us for comparison.}. 
%For the CPLEX solver, we use its native Python API of CPLEX v22.1.1 with Python 3.10.16 (as the API does not yet support Python 3.11), while we use Google OR-Tools CP-SAT with Python 3.11.11. %, which allowed for a comparison with an open-source solver. 
For ILP, we used the CPLEX v22.1.1 Python API with Python v3.10.16,
%(as the API does not yet support Python v3.11)
while we used Google OR-Tools CP-SAT with Python v3.11.11.
For CPLEX and Google OR-Tools, we used the default settings. 

The ILP algorithms operate on two modalities:
\begin{enumerate*}[label=\alph*)]
    \item daily schedule: scheduling is performed at the end of each day, taking into account all patients who arrived that day. Then, the machines' capacity is updated, and the next day is scheduled until all patients have been treated. %On the other end, 
    \item all patients' schedules at a time: we take all the patients from the dataset and schedule them all at once to find what the best schedule would be if we knew in advance all the patients' arrivals.
    %a more detailed description can be found in the Methodology section.
\end{enumerate*}

Due to time constraints, the CPLEX and CP-SAT experiments were run only once, with the first seed of the series used for the SA (see later). We considered two different time limits for these exact methods: 1) 1 hour (the same time limit used by \cite{frimodig2023comparing}) for both modalities; and 2) 200 seconds (a value closer to the average time needed by the SA algorithm--see later) for the first modality only, since the second would not have led to feasible solutions. 
%The two configurations--the all-patient schedule with a 3600-second limit and the daily batch configuration with a 200-second limit--tested using the CP-SAT solver produced poor results for the former and no results for the latter. Given these outcomes, we decided not to report or discuss them further.

The heuristic algorithms are deterministic: given the same daily set of patients, they always apply the same priority-based ordering and the same scheduling rules. Therefore, running them multiple times is unnecessary. These heuristics operate in an online modality: not only do they process patients day-by-day, but within each day, they schedule patients sequentially, making irrevocable decisions one patient at a time.
%The heuristic algorithms are deterministic: they operate on an online modality by processing patients day by day, and within each day, scheduling patients sequentially, making irrevocable decisions one patient at a time. Thus, given the same daily set of patients, they always apply the same priority-based ordering and the same scheduling rules. Therefore, running them multiple times is unnecessary. 
In contrast, each of the SA implementations (\textit{plain}, \textit{daily}, and \textit{reheat}) has been run five times with five different seeds--i.e., 198743, 3947394, 50343784, 93790244, 234720309. 
%both for the daily schedule and for all patients' schedules at a time, with both \textit{plain} and \textit{reheat} implementations. 
The parameters of each SA configuration have been tuned with iRace 4.3 \cite{lopez2016irace} on a subset of instances. 
%The specific values can be found in the Supplementary Material.

If any of the algorithms cannot find a solution within the considered time limit, the Obj. Sum for the unfeasible instance is set equal to the maximum Obj. Sum found over all instances, to allow for a consistent scale of the Obj. Sum values for each algorithm.

For all formulations, the number of days (bins for the heuristics) considered in the time horizon for patient scheduling is determined by each problem instance.

\paragraph{Statistical evaluation}
Following \cite{demvsar2006statistical} guidelines, to statistically compare the different methodologies, we used the Friedman test to reject the null hypothesis (that the methodologies compared are not statistically different), and we performed a pairwise analysis using the Wilcoxon signed-rank test. For overall ranking, we used Critical Difference (CD) diagrams \cite{IsmailFawaz2018deep}, which show the rank from the worst to the best--left to right--of each methodology compared, and draw a horizontal tick line between methodologies that are not significantly different in terms of the metric compared.
%For each comparison, all different instances and combinations of weighted Obj. Sum have been tested for differences between methodologies.

\subsection{Results}

\begin{figure*}[t!]
    \centering
    \includegraphics[width=.92\linewidth]{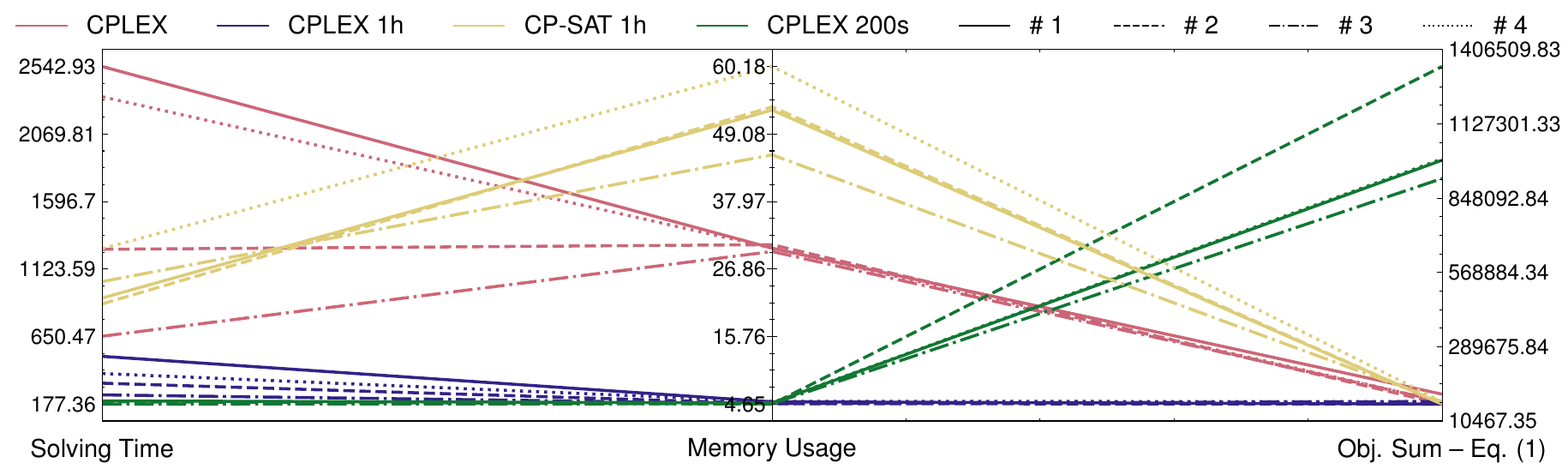}
%    \vspace{-0.2cm}
    \caption{
    Parallel coordinate plot of solving time (seconds), memory usage (GB), and best Obj. Sum value (see \Cref{eq:objective_sum}), for each weighted objective combination (\#1–\#4), averaged over all instances. Colours represent solvers; line styles, weighted combinations.
    %Parallel coordinate plot showing Solving time (seconds) on the first vertical axis, Memory usage (GB) on the second vertical axis, and best value found for Obj. Sum (see \Cref{eq:objective_sum}) on the third vertical axis, for each weighted combination of objectives (from \#1 to \#4), averaged across all instances. Each colour corresponds to a different solver. Each line style corresponds to the different weighted combinations considered.
    }
    %ILP parallel plot comparing performance of exact methods implementation on all evaluation metrics
    \label{fig:ILP-pp}
\end{figure*}

\paragraph{ILP} 
\Cref{fig:ILP-pp} presents the performance comparison of the ILP implementations. \textit{CPLEX} refers to the all-patients schedule modality with a 1-hour time limit, while all the others will refer to the daily batch modality; the time limit is indicated in the legend, i.e., \textit{CPLEX 1h} uses a time limit of 1 hour. We remark that ILP CP-SAT, for all patients' schedules with a 1-hour time limit, and for the daily schedule with a 200-second time limit, produced poor-quality solutions for the former and no solution at all for the latter; therefore, we decided not to report them in the plots or discuss them further.

From the figure, it can be seen that the \textit{CPLEX 1h} configuration is the one that can find better results in terms of solution quality and memory usage metrics, while \textit{CPLEX 200s} is obviously--by construction--the best in terms of solving time, but with really poor results in terms of solution quality. The other two configurations use both more solving time and memory, but yield solutions of better quality. Overall, across all the weighted objective combinations and instances--encompassing a total of 320 experiments--\textit{CPLEX} found a solution for 250 experiments, \textit{CP-SAT 1h} for 293 experiments, \textit{CPLEX 200s} for 112, while \textit{CPLEX 1h} solved every instance for every weighted objective combination. 
For subsequent comparisons, only \textit{CPLEX 1h} and \textit{CPLEX 200s} have been chosen. The first is the best one on almost all the evaluation metrics, while the latter is the one with a solving time comparable to the SA.

\paragraph{Heuristics}
The heuristics are the algorithms with better solving time and memory usage. RTSP FF can solve the instances over all weighted Obj. Sum combinations in 0.5 seconds, using 0.1 GB of memory on average, while RTSP BF takes 3.3 seconds and uses 0.1 GB of memory on average. Since they do not actively search for an optimal solution, they can reach sub-optimal values, as shown in \Cref{tab:results-summary}, where the Obj. Sum value, averaged across all instances, is presented on the first row with each combination on each sub-row.

\begin{comment}
\begin{table}[t!]
\centering
{\scriptsize
\begin{tabular}{ccc}
\toprule
\textbf{Comb.} & \textbf{RTSP FF} & \textbf{RTSP BF} \\
\midrule
\#1 & 14174.11 $\pm$ 4372.46 & 25043.79 $\pm$ 4428.29 \\
\#2 & 13227.93 $\pm$ 4408.55 & 14866.04 $\pm$ 4410.22 \\
\#3 & 23199.61 $\pm$ 8165.05  & 25714.41 $\pm$ 7904.76 \\
\#4 & 32060.68 $\pm$ 9474.60 & 48277.54 $\pm$ 10466.17 \\
\bottomrule
\end{tabular}
}
\caption{Average $\pm$ std. dev. Obj. Sum values across all instances for RTSP FF and RTSP BF.}
\label{tab:heuristics-obj}
\end{table}

\begin{table*}[h!]
\centering
\caption{Average $\pm$ std. dev. Obj. Sum values across all instances for ILP methods, RTSP FF, and RTSP BF.}
{\scriptsize
\begin{tabular}{ccccc}
\toprule
\textbf{Comb.} & \textbf{CPLEX 1h} & \textbf{CPLEX 200s} & \textbf{RTSP FF} & \textbf{RTSP BF} \\
\midrule
\#1 & \textbf{10530.67} $\pm$ 2987.74 & 1019050.35 $\pm$ 177779.48 & 14174.11 $\pm$ 4372.46 & 25043.79 $\pm$ 4428.29 \\
\#2 & \textbf{10563.88} $\pm$ 3011.81 & 1406509.83 $\pm$ 595054.10 & 13227.93 $\pm$ 4408.55 & 14866.04 $\pm$ 4410.22 \\
\#3 & \textbf{20958.55} $\pm$ 5910.51 & 943360.09 $\pm$ 383715.77 & 23199.61 $\pm$ 8165.05  & 25714.41 $\pm$ 7904.76 \\
\#4 & \textbf{21163.32} $\pm$ 5994.59 & 1024139.58 $\pm$ 188653.08 & 32060.68 $\pm$ 9474.60 & 48277.54 $\pm$ 10466.17 \\
\bottomrule
\end{tabular}
}
\label{tab:heuristics-obj}
\end{table*}
\end{comment}

\paragraph{Simulated Annealing}
The SA \textit{plain}, \textit{daily}, and \textit{reheat} implementations have been tested on all different instances and combinations of weighted Obj. Sum. 
%In order to have a statistical comparison of the SA implementation between the heuristics, the ILP methods and themselves, we used Critical Difference (CD) diagrams \cite{IsmailFawaz2018deep}. 
The results of these experiments are reported in \Cref{fig:all-sa} (a), where it is clear that the \textit{plain} implementation is not only statistically different from the other for both starting solutions, but it is also the best one in terms of solution quality. Therefore, for the subsequent comparisons, we focus only on the \textit{plain} SA.

\subsection{Critical discussion of the results}
\begin{figure*}[ht!]
    \centering
    \addtolength{\tabcolsep}{-0.45em}
    \begin{tabular}{>{\tiny}c@{}c@{}>{\tiny}c@{}c@{}}
        \textbf{a)} & %\begin{subfigure}{0.95\columnwidth}
    \includegraphics[align=c,width=0.95\columnwidth]{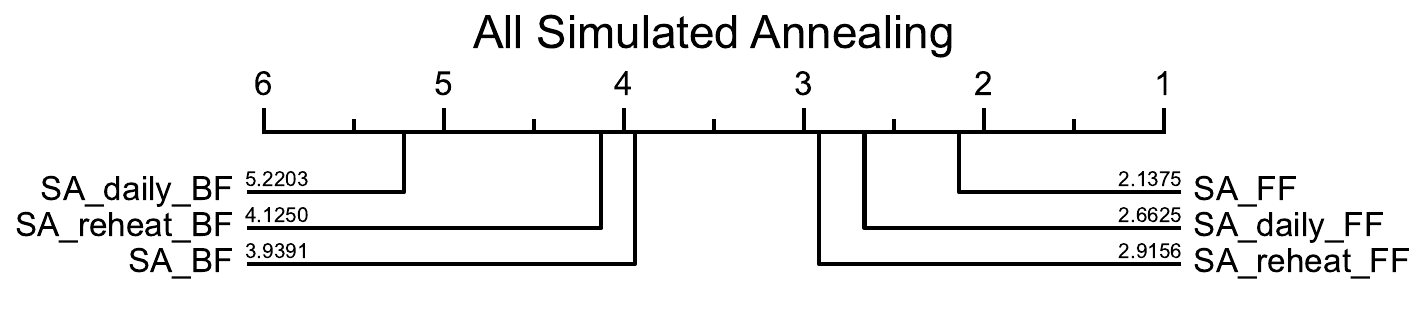}
    %\end{subfigure} 
    &
        \textbf{d)} & %\begin{subfigure}{0.95\columnwidth}
    \includegraphics[align=c,width=0.95\columnwidth]{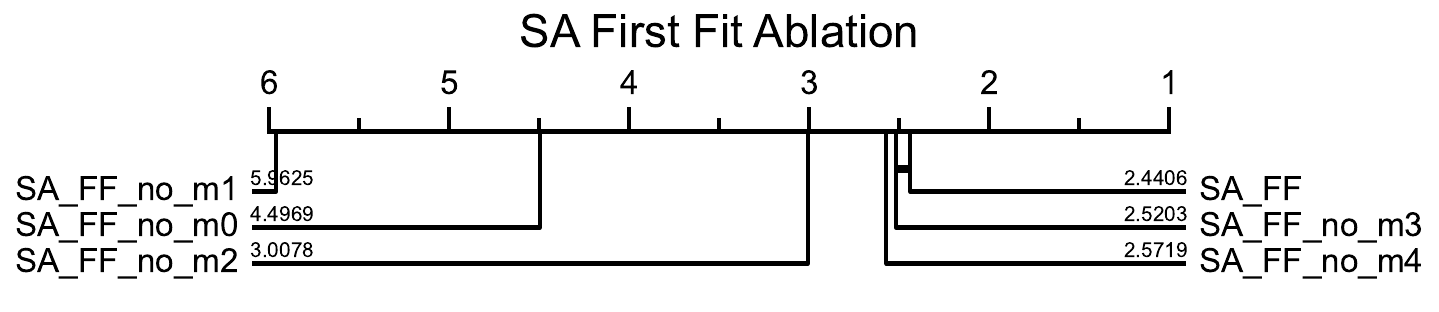}
    %\end{subfigure}
    \\
        \textbf{b)} & %\begin{subfigure}{0.95\columnwidth}
    \includegraphics[align=c,width=0.95\columnwidth]{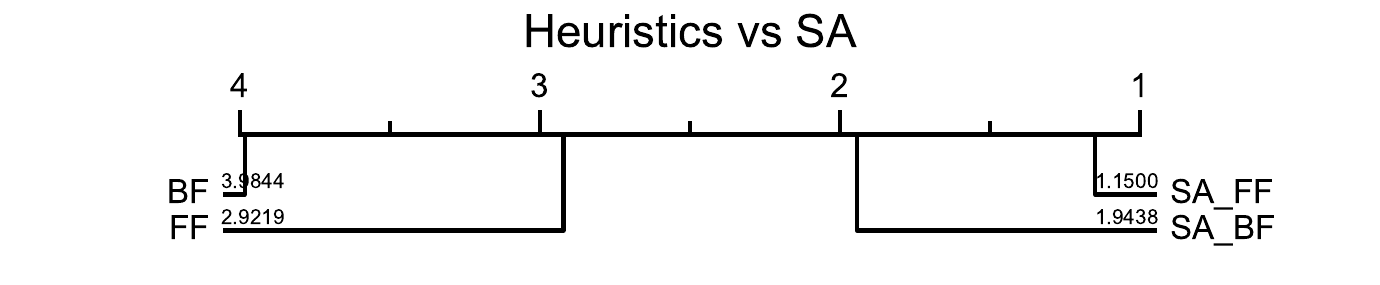}
    %\end{subfigure} 
    &
        \textbf{e)} & %\begin{subfigure}{0.95\columnwidth}
    \includegraphics[align=c,width=0.95\columnwidth]{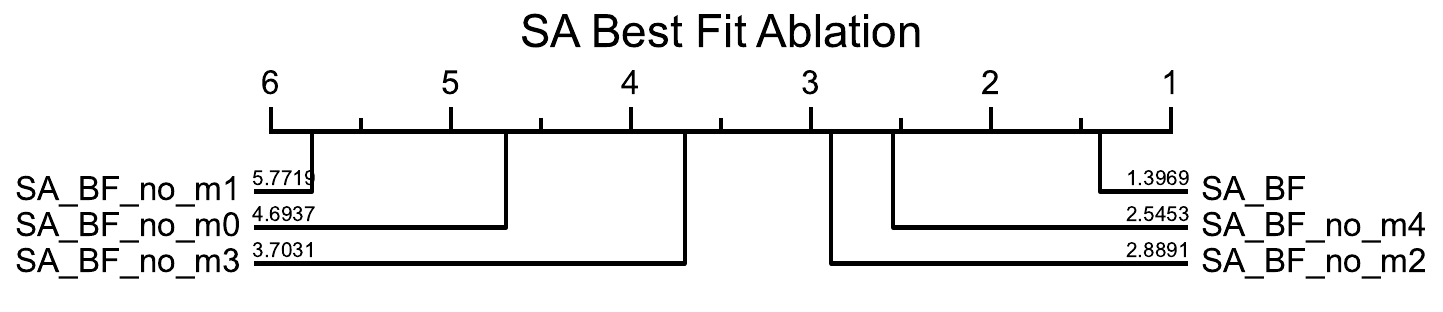}
    %\end{subfigure}
    \\
        \textbf{c)} & %\begin{subfigure}{0.95\columnwidth}
    \includegraphics[align=c,width=0.95\columnwidth]{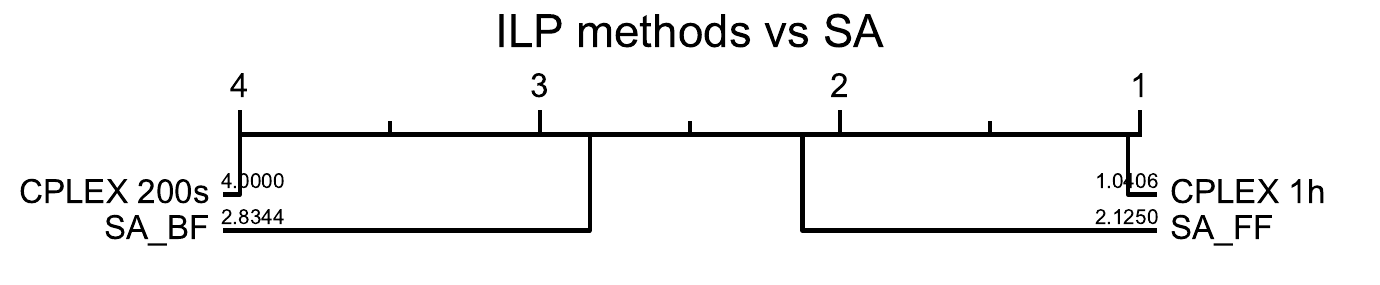}
    &
    %\end{subfigure}
    \end{tabular}
    \caption{CD diagrams for a) all SA implementations; b) heuristics vs. SA; c) SA vs. ILP methods; d) SA starting from FF heuristic and e) SA starting from BF heuristic.}% ablations
    \label{fig:all-sa}
\end{figure*}

\Cref{fig:all-sa} (b) shows that the \textit{plain} implementation of the SA effectively improves the quality of the solution obtained from the heuristic algorithms. Moreover, as RTSP FF yields better solutions than RTSP BF, so does the SA starting from the former (SA\_FF) compared to the SA starting from the latter (SA\_BF).

The CD diagram in \Cref{fig:all-sa} (c) compares the \textit{plain} SA and the two ILP methods chosen for comparison, i.e., \textit{CPLEX 1h} and \textit{CPLEX 200s}. This diagram clearly shows that both the SA \textit{plain} algorithms, i.e., the one starting from RTSP FF solutions and the one starting from RTSP BF solutions, can outperform the CPLEX solver when using a comparable amount of time, i.e., the \textit{CPLEX 200s} setting. On the other hand, when CPLEX is given more time to find the solution, i.e., \textit{CPLEX 1h}, it can outperform SA algorithms in terms of solution quality. \Cref{tab:results-summary} presents a more detailed comparison of the evaluation metrics results for \textit{CPLEX 1h}, \textit{CPLEX 200s}, \textit{SA\_FF}, and \textit{SA\_BF}.

\begin{table*}[ht!] 
\centering 
\caption{Comparison of SA and ILP performance across Obj. Sum, time, memory, and heuristic baselines. Values are presented as mean $\pm$ std. dev. across all instances and across seeds for the SA cases. The Obj. Sum. values must be multiplied by $10^{3}$.} {\scalebox{0.99}{
\renewcommand{\arraystretch}{0.99} 
%\addtolength{\tabcolsep}{-0.48em}
%\begin{tabular}{l@{\ \ }lcccccc} 
% \begin{tabular}{l@{\ \ }lr@{ $\pm$ }lr@{ $\pm$ }lr@{ $\pm$ }lr@{ $\pm$ }lr@{ $\pm$ }lr@{ $\pm$ }l} 
\begin{tabular}{
l l
r@{ $\pm$ }l @{\hspace{0.95em}}
r@{ $\pm$ }l @{\hspace{0.95em}}
r@{ $\pm$ }l @{\hspace{0.95em}}
r@{ $\pm$ }l @{\hspace{0.95em}}
r@{ $\pm$ }l @{\hspace{0.95em}}
r@{ $\pm$ }l
}
\toprule 
\multicolumn{2}{l}{\textbf{Comb.}} & \multicolumn{2}{c}{\textbf{CPLEX 1h}} & \multicolumn{2}{c}{\textbf{CPLEX 200s}} & \multicolumn{2}{c}{\textbf{RTSP FF}} & \multicolumn{2}{c}{\textbf{RTSP BF}} & \multicolumn{2}{c}{\textbf{SA\_FF}} & \multicolumn{2}{c}{\textbf{SA\_BF}} \\ 
\midrule 

%===================== OBJECTIVE ===================== 
\parbox[t]{2mm}{\multirow{4}{*}{\rotatebox[origin=c]{90}{\scalebox{0.8}{\textbf{Obj. Sum}}}}}
& \#1 & \textbf{10.53} & 2.99 & 1019 & 178 & 14.17 & 4.37 & 25.04 & 4.43 & 12.67 & 3.93 & 13.93 & 4.04 \\ 
& \#2 & \textbf{10.56} & 3.01 & 1407 & 595 & 13.23 & 4.41 & 14.87 & 4.41 & 11.28 & 3.08 & 11.38 & 3.18 \\ 
& \#3 & \textbf{20.96} & 5.91 & 943 & 384 & 23.20 & 8.17 & 25.71 & 7.90 & 21.94 & 6.21 & 22.55 & 7.34 \\ 
& \#4 & \textbf{21.16} & 5.99 & 1024 & 189 & 32.06 & 9.47 & 48.28 & 10.47 & 26.58 & 8.10 & 27.30 & 7.66 \\ \midrule 

% ===================== TIME ===================== 
\parbox[t]{2mm}{\multirow{4}{*}{\rotatebox[origin=c]{90}{\scalebox{0.8}{\textbf{Time}}}}}
&  \#1 & 512.12 & 482.08 & 199.76 & \ 2.16 & \textbf{0.53} & 0.12 & 3.34 & 0.001 & 24.91 & 0.98 & 40.21 & 3.49 \\
&  \#2 & 324.05 & 285.04 & 193.99 & \ 9.80 & \textbf{0.53} & 0.16 & 3.29 & 0.001 & 32.49 & 3.28 & 33.34 & 2.06 \\
&  \#3 & 241.98 & \ 97.73 & 177.36 & 16.63 & \textbf{0.53} & 0.13 & 3.30 & 0.001 & 23.43 & 2.40 & 29.20 & 1.69 \\
&  \#4 & 391.73 & 321.68 & 195.60 & \ 9.79 & \textbf{0.54} & 0.16 & 3.26 & 0.001 & 26.78 & 2.50 & 36.99 & 2.94 \\ 
\midrule 

% ===================== MEMORY ===================== 
\parbox[t]{2mm}{\multirow{4}{*}{\rotatebox[origin=c]{90}{\scalebox{0.8}{\textbf{Memory}}}}}
& \#1 & 5.03 & 1.20 & 4.84 & 1.12 & 0.099 & 0.40 & 0.099 & 0.001 & \textbf{0.04} & 0.0004 & \textbf{0.04} & 0.0004 \\
& \#2 & 4.73 & 1.17 & 4.65 & 1.05 & 0.099 & 0.40 & 0.099 & 0.001 & \textbf{0.05} & 0.0021 & \textbf{0.05} & 0.0019  \\
& \#3 & 4.74 & 1.17 & 4.71 & 1.15 & 0.099 & 0.40 & 0.099 & 0.001 & \textbf{0.05} & 0.0028 & \textbf{0.05} & 0.0011 \\
& \#4 & 5.09 & 1.41 & 4.89 & 1.05 & 0.099 & 0.40 & 0.099 & 0.001 & \textbf{0.04} & 0.0004 & \textbf{0.04} & 0.0004 \\ 
\bottomrule 
\end{tabular} 
}
}
\label{tab:results-summary} 
\end{table*}

\paragraph{Ablation}
To better evaluate SA performance, an ablation study has been performed. The study is conducted by analysing how the SA performs when using all neighbourhoods except one \cite{fawcett2016analysing}. To exclude a move, the probability associated with it is zeroed, while the remaining are scaled to sum to 1, maintaining their in-between ratio. 
A Friedman test with a significance level $\alpha = 0.05$ was performed against the plain SA, and all ablation experiments were performed with the starting solution obtained by both RTSP FF and RTSP BF, to determine if they perform equally. For the SA\_FF experiments, the \textit{p-value} is $1.60e^{-194}$, while for the SA\_BF the \textit{p-value} is $4.33e^{-222}$, rejecting in both cases the null hypothesis and therefore stating that the two sets are statistically different.

\begin{comment}
\begin{figure}[ht!]
    \centering
    \addtolength{\tabcolsep}{-0.45em}
    \begin{tabular}{>{\small}c@{}c@{}}
        \textbf{a)} & %\begin{subfigure}{0.95\columnwidth}
    \includegraphics[align=c,width=0.95\columnwidth]{Images/CriticalDifference/cd-diagram_all_SA.pdf}
    %\end{subfigure} 
    \\
        \textbf{b)} & %\begin{subfigure}{0.95\columnwidth}
    \includegraphics[align=c,width=1\columnwidth]{Images/CriticalDifference/cd-diagram_he_sa.pdf}
    %\end{subfigure} 
    \\
        \textbf{c)} & %\begin{subfigure}{0.95\columnwidth}
    \includegraphics[align=c,width=1\columnwidth]{Images/CriticalDifference/cd-diagram_sa_ilp.pdf}
    %\end{subfigure}
    \end{tabular}
    \caption{CD diagrams (from the top: all SA configurations; heuristics vs. SA; SA vs. ILP methods).}
    \label{fig:all-sa}
\end{figure}

\begin{figure}[ht!]
    \centering
\begin{tabular}{@{}c}
%    \begin{subfigure}{0.95\columnwidth}
        \includegraphics[width=0.95\columnwidth]{Images/CriticalDifference/cd-diagram_ff_ablation.pdf}
%    \end{subfigure}
\\
%    \begin{subfigure}{0.95\columnwidth}
        \includegraphics[width=0.95\columnwidth]{Images/CriticalDifference/cd-diagram_bf_ablation.pdf}
%    \end{subfigure}
\end{tabular}
    \caption{CD diagrams (SA starting from the heuristics).}
    \label{fig:sa-ablation}
\end{figure}
\end{comment}

Following the Friedman test, we further analyse the configurations using the CD diagram in \Cref{fig:all-sa} (d) and (e), 
comparing the ablations of SA starting from the RTSP FF (\textit{SA\_FF}) and RTSP BF (\textit{SA\_BF}) solutions.
For each ablation, the SA variant in which the probability of performing a given move is set to zero is denoted as \textit{no\_m\#}, where \# refers to the move index defined in the description of the SA methodology.

\begin{table*}[ht!]
\centering
\caption{
Percentage of accepted and improving moves across all weighted Obj. Sum combinations and SA\_FF implementations. P, D, and R stand for, respectively, Plain, Daily, and Reheat SA implementations.}
{%\scriptsize
\scalebox{0.99}{
\setlength{\tabcolsep}{3pt}
\renewcommand{\arraystretch}{0.99}
\begin{tabular}{l@{\ \ }lcccccccccccc}
\toprule
&& \multicolumn{3}{c}{\textbf{Comb. \#1}} 
& \multicolumn{3}{c}{\textbf{Comb. \#2}}
& \multicolumn{3}{c}{\textbf{Comb. \#3}}
& \multicolumn{3}{c}{\textbf{Comb. \#4}} \\
\cmidrule(lr){3-5}\cmidrule(lr){6-8}\cmidrule(lr){9-11}\cmidrule(lr){12-14}

 && P & D & R & P & D & R & P & D & R & P & D & R \\
\midrule

% Obj. Sum & 12679.62 & 12633.61 & 12839.66
%      & 11273.12 & 10936.36 & 11487.69
%      & 21881.19 & 21499.07 & 22044.45
%      & 26568.50 & 27325.34 & 26782.55 \\

% Time & 26.31 & 55.26 & 29.76
%      & 34.39 & 85.15 & 29.41
%      & 25.07 & 68.82 & 20.20
%      & 28.18 & 48.01 & 26.54 \\

% Memory & 0.04 & 0.08 & 0.09
%        & 0.05 & 0.08 & 0.09
%        & 0.05 & 0.08 & 0.09
%        & 0.04 & 0.08 & 0.09 \\
% \midrule

\parbox[t]{2mm}{\multirow{4}{*}{\rotatebox[origin=c]{90}{\scalebox{0.8}{\textbf{Accepted  }}}}}
& m0
& 6.81 & 7.17 & 6.20
& 27.33 & 1.96 & 2.92
& 4.42 & 3.75 & 4.81
& 27.27 & 28.15 & 27.23 \\

& m1
& 2.51 & 2.78 & 3.21
& 83.00 & 79.17 & 79.45
& 52.43 & 52.19 & 55.80
& 2.90 & 2.98 & 3.33 \\

& m2
& 2.12 & 5.12 & 1.73
& 1.69 & 3.91 & 1.71
& 0.98 & 3.08 & 1.18
& 2.96 & 6.34 & 2.83 \\

& m3
& 5.63 & 14.11 & 5.62
& 91.34 & 94.05 & 91.97
& 66.75 & 69.51 & 67.02
& 5.14 & 10.17 & 5.12 \\

& m4
& 0.03 & 1.03 & 0.05
& 0.02 & 1.05 & 0.01
& 0.02 & 1.06 & 0.03
& 0.04 & 1.06 & 0.05 \\

\midrule

\parbox[t]{2mm}{\multirow{4}{*}{\rotatebox[origin=c]{90}{\scalebox{0.8}{\textbf{Improving}}}}}
& m0
& 1.03 & 0.99 & 0.86
& 0.90 & 0.49 & 1.09
& 0.40 & 0.40 & 0.47
& 1.08 & 1.63 & 1.06 \\

& m1
& 0.12 & 0.19 & 0.15
& 0.01 & 0.00 & 0.001
& 0.02 & 0.01 & 0.03
& 0.16 & 0.25 & 0.18 \\

& m2
& 0.37 & 0.46 & 0.27
& 0.47 & 0.19 & 0.49
& 0.14 & 0.12 & 0.16
& 1.35 & 1.64 & 1.28 \\

& m3
& 0.07 & 0.10 & 0.05
& 0.003 & 0.00 & 0.00
& 0.02 & 0.01 & 0.02
& 0.09 & 0.14 & 0.06 \\

& m4
& 0.03 & 0.002 & 0.05
& 0.02 & 0.003 & 0.01
& 0.02 & 0.004 & 0.03
& 0.04 & 0.009 & 0.05 \\

\bottomrule
\end{tabular}
}
}
\label{tab:move-perc}
\end{table*}

For \textit{SA\_FF}, the only move whose removal does not produce results statistically different from the full algorithm is \textit{Swap Machines (m3)}. This is likely due to the fact that many patients can be treated on multiple beam-matched machines, making a large portion of machine swaps ineffective. When two patients exchange two beam-matched machines that are also among their preferred options, the weighted objective value remains unchanged (see the $f_5$ and $f_6$ objectives in the Problem Formulation section).
As expected, disabling any of the other moves leads to statistically worse outcomes.

For \textit{SA\_BF}, the CD diagram shows that none of the ablations are statistically equivalent to the whole algorithm: removing a single move consistently leads to inferior performance.

\paragraph{Statistics}
To understand how SA explores the neighbourhood, we analyse the percentages of accepted and improving moves for each weighted sum combination and SA variant, using First Fit as starting heuristic since it yields better performance. Accepted moves are those whose resulting solutions are retained by SA, while improving moves are those that actually produce better solutions.
%To better understand how SA explores the neighbourhood, the percentages of accepted and improving moves are analysed for each weighted sum combination and for each variant of SA--with FF as start heuristic since is the one performing better. Accepted moves indicate the frequency with which the move leads to a solution that is selected and retained by SA, while improving moves indicate the frequency with which that move actually produces a better solution.

As shown in \Cref{tab:move-perc}, \textit{m0} is more relevant for combinations \#2--only for SA \textit{plain} implementation--and \#4, with 27-28\% accepted moves, suggesting that the move contributes more to search diversification. On the other hand, despite being the most improved move among all the others, the improved percentage is quite low--around 0.5-1.5\%--suggesting that, in general, the First Fit heuristic is already finding near-to-optimal solutions and, in particular, that \textit{m0} is the most effective one. \textit{m1} is the predominant move for combinations \#2 and \#3, with very high accepted moves' percentages (up to 83\%) but low improving moves' percentages, meaning that it supports exploration over direct improvement. \textit{m2} is slightly higher for combinations \#1 and \#4 for the SA \textit{daily} implementation, while it remains at lower values for the others. It's improving percentages are the second best, after \textit{m0}, especially for combination \#4, suggesting the move tends to improve more than the others when accepted. 
\textit{m3} shows a large discrepancy between the percentages of accepted and improving moves, especially for combinations \#2 and \#3.
This indicates that many machine swaps move patients among beam-matched machines without improving the quality of the solution.
\textit{m4} shows very low percentages in terms of both accepted and improving moves, indicating that the move is rarely selected and provides little search value.

\paragraph{Overall considerations}
The experimental evaluation highlights clear trade-offs between exact methods and heuristic/metaheuristic approaches for solving the RTSP.
First, ILP can reach optimal solutions when enough runtime is available, although it requires more memory. However, its performance drops rapidly under strict time limits.
Second, the proposed RTSP heuristics--despite not performing an active search--are able to produce good solutions with very low computational cost.
Third, the SA algorithm can further improve the heuristic solutions while keeping runtime short and memory usage low. The analysis also shows the importance of the current neighbourhoods, and there is a potential for new neighbourhoods to be explored.
In general, heuristics and SA achieve a good balance between solution quality and computational effort. Although ILP remains the most reliable choice for optimality when runtime is not restricted, SA offers a valid alternative, producing near-optimal schedules with a significant reduction in both runtime and memory usage.

\section{Conclusion and Future Work}
\label{sec:conclusions}
The RTSP has been widely tackled with exact methods and metaheuristics. However, the existing literature reveals a lack of work focused on simpler heuristics. In this paper, we developed two novel greedy heuristics for the RTSP, which are then combined with a custom SA algorithm. We compared the heuristics and three SA implementations with different exact methods. The experimental results on a public dataset simulating a real-world arrival flow to a large cancer centre demonstrate that simple heuristics such as RTSP FF and RTSP BF can drastically reduce both computational time and memory usage. Moreover, the SA starting from the solutions found by heuristics can further improve the solution quality while keeping the computational time and memory usage drastically low compared to the exact methods. This makes our approach attractive for hospitals where physicians may have access to limited computational resources and cannot allocate budget for the purchase of expensive solver licences. However, this computational efficiency comes at the cost of achieving suboptimal solutions for the four weight combinations investigated in comparison with the exact methods run for a longer time. To address this limitation, future work will focus on refining the SA moves to increase the improvement percentage, aiming at neighbourhoods that not only explore but also exploit the solution quality. Moreover, we will investigate new heuristics, extending the RTSP BF and FF to better take into account the objective function.

\bibliographystyle{IEEEtran}
\bibliography{bibliography}
\end{document}